\documentclass[11pt]{report}

\usepackage[T2A]{fontenc}
\usepackage[cp1251]{inputenc}
\usepackage{amsmath}
\usepackage{amssymb,latexsym}
\usepackage{oldlfont}
\usepackage{amsfonts}

\begin{document}

\begin{center}

{{\Large \bf 
On extension of partial orders to total preorders\\  with prescribed
symmetric part}\footnote{The research was supported by the National
Program of Fundamental Researches of Belarus under grant
``Convergence -- 1.4.03''. }}

\vspace{7mm}

{\large \bf Dmitry V. Akopian and Valentin V. Gorokhovik}
\vspace{0.3cm}

\vspace{2mm}

{\small Institute of Mathematics,\\ The National Academy of Sciences of Belarus, \\
Surganova st., 11, Minsk 220072, Belarus\\ e-mail: gorokh@im.bas-net.by}

\end{center}

\vspace{7mm}

{{\bf Abstract} For a partial order $\preceq$ on a set $X$ and an equivalency relation
$S$ defined on the same set $X$ we derive a necessary and sufficient condition for the
existence of such a total preorder on $X$ whose asymmetric part contains the asymmetric
part of the given partial order $\preceq$ and whose symmetric part coincides with the
given equivalence relation $S.$ This result generalizes the classical Szpilrajn theorem
on extension of a partial order  to a perfect (linear) order. }

\vspace{3mm}

{{\bf Key words} partial order, preorder, extension, Szpilrajn theorem,
Dushnik-Miller theorem.}

\vspace{3mm}

{{\bf MSC 2010}\, Primary: 06A06 Secondary: 06A05, 91B08}

\vspace{7mm}

\hspace*{3pt} Let $X$ be an arbitrary nonempty set and let $G \subset X\times X$ be a
binary relation on $X.$ A binary relation $ G \subset X \times X$ is called {\it a
preorder} if it is reflexive ($(x,\,x) \in G\,\,\forall\,\,x \in X$) and transi\-tive
($(x,\,y) \in G,\,(y,\,z) \in G \Rightarrow (x,\,z) \in G \,\,\forall\,\,x,y,z \in X$).
If in addition a preorder $G \in X \times X$ is antisymmetric ($(x,\,y) \in G,\,(y,\,x)
\in G \Rightarrow x = y\,\,\forall\,\,x,y \in X$), then it is called {\it a partial
order}. A total partial order is called {\it a perfect} (or {\it linear}) {\it order}.
(A binary relation $G \subset X\times X$ is {\it total} if for any $x,\,y \in X$ either
$(x,\,y) \in G$ or $(y,\,x) \in G$ holds.) In the sequel, a preorder will be preferably
denoted by the symbol $\precsim$ whereas a partial order as well as a perfect order by
the symbol $\preceq.$

\vspace{2mm}

One of the key results of the theory of ordered sets is the following theorem proved by
E. Szpilrajn in 1930 [1].

\vspace{3mm}

{\bf Тheorem 1} (E. Szpilrajn [1])\,\, {\it For every partial order $\preceq \subset X
\times X$ there exists a perfect extension, i.~e., there exists a perfect order
$\preceq' \subset X \times X$ such that $\preceq \subset \preceq'.$ Moreover, for any
pair of elements $a,\,b \in X$ such that $(a,\,b)\not \in \preceq$ and $(b,\,a) \not \in
\preceq$ a perfect extension $\preceq' \subset X \times X$ for the partial order
$\preceq$ can be chosen in such a way that $(a,\,b) \in \preceq'.$}

\vspace{3mm}

In 1941 E. Dushnik and B. Miller proved the following strengthening of the Szpilrajn
theorem.

\vspace{3mm}

{\bf Theorem 2} (E. Dushnik, B. Miller [2])\,\, {\it Every partial order $\preceq \subset X
\times X$ is the intersection of all its perfect extensions.}

In the recent literature the Szpilrajn theorem and the Dushnik--Miller theorem and their
proofs can be found in the monographs [3,4]. The generalizations of the Szpilrajn
theorem to the case when partial orders and perfect orders extending them are defined on
groups, rings and some other algebraic systems and are compatible with their algebraic
operations are presented in the monograph of L. Fuchs [5].
Due to the duality between compatible perfect orders defined on a real vector space $X$
and semispaces of $X$ (the cones of positive elements of compatible perfect orders are
complements of semispaces at zero) it follows from the results of V.~Klee devoted to
semispaces [6] that any compatible partial order defined on a real vector space $X$ can be extended to a compatible perfect order. For relations defined  on topological spaces the conditions under which there exist continuous total preorders extending partial orders were obtained by G. Bosi
and G. Herden [7,\,8]. The results of the studies devoted to the existence of utility
functions for partial orders (see [3,\,9] as well the monographs [10,\,11] and bibliography cited
there) can also be considered as generalizations of the Szpilrajn theorem.

Every binary relation $G$ on $X$ can be presented as the disjoint union $G = P_G\cup S_G$ ($P_G\cap S_G = \varnothing$) of its asymmetric part $P_G:= \{(x,\,y)
\in G\,|\,(y,\,x) \not \in G\}$ and its symmetric part $S_G:=\{(x,\,y) \in G\,|\,(y,\,x)
\in G\}.$
If $G$ is a preorder then its symmetric part $S_G$ is reflexive, symmetric and
transitive and, consequently, in that case $S_G$ is an equivalency relation on $X,$
which is reduced to the equality relation when $G$ is  a partial order. The asymmetric
part of a preorder (and, in particular, the asymmetric part of a partial order) is an
asymmetric and transitive binary relation. On the other hand, the union of any
asymmetric and transitive binary relation with the equality relation is a partial order.
Thus, there exists the one-to-one correspondence between partial orders and asymmetric
and transitive binary relations. Note that different preorders can have the same
asymmetric part.

Let $\preceq$ be a partial order on $X$ and $S$ an equivalency relation defined on the same set $X.$

A total preorder $\precsim \subset X \times X$ will be referred to as  {\it a total
preorder $S-$extension of the partial order $\preceq$} if the asymmetric part of
$\precsim$ contains the asymmetric part of the given partial order $\preceq$ and the
symmetric part of $\precsim$ coincides with the equivalency relation $S,$ that is, if
$P_\preceq \subset P_{\precsim}$ and $S_{\precsim} = S.$

The main purpose of this paper is to derive for a given partial order $\preceq$ and a
given equivalency relation $S$ a necessary and sufficient condition for the existence of
a total preorder $S-$extension of $\preceq.$

In the case when $S$ is the equality relation on $X,$ i.~e., when $S = E :=\{(x,\,y) \in X \times X\,|\, x = y\},$ due to the Szpilrajn theorem, such an
extension exists for any partial order $\preceq.$ As it will be shown below in the
general case the required extension exists if and only if the partial order $\preceq$ and
the equivalency relation $S$ are compatible in some way. Thus the main results of the
paper can be considered as a generalization of  the Szpilrajn theorem.

Let us begin with consideration of a particular case. Assume that a partial order
$\preceq$ and an equivalency relation $S$ hold the additional condition
\begin{equation}\label{e1}
P_\preceq \circ S = S \circ P_\preceq = P_\preceq
\end{equation}
(the symbol  $\circ$ denotes the composition of binary relations).

It immediately follows from \eqref{e1} that the union $\preceq \cup S$ is a preorder on
$X$ the symmetric part of which coincides with $S.$ Let $X/S$ be the quotient of $X$
with respect to the equivalency relation $S$ and let $T$ be the quotient of the preorder
$\preceq \cup S$ with respect to $S.$ Since $T$ is a partial order on $X/S,$ due to the
Szpilrajn theorem, $T$ can be extended to a perfect order $Q$ on $X/S.$ Setting $x
\precsim y \Longleftrightarrow [x]_S Q [y]_S$ (here $[x]_S$ and $[y]_S$ stands for the
equivalency classes of $S$ containing $x$ and $y,$ respectively) we obtain the total
preorder $\precsim$ on $X$ which is a total preorder $S-$extension of $\preceq.$

Thus the following theorem generalizing both the Szpilrajn theorem and the
Dushnik--Miller theorem is true.

{\bf Theorem 3}\,\, {\it Let $\preceq$ be a partial order on $X.$ For any equivalency
relation $S$ on $X$ which satisfies condition \eqref{e1}, there exists a total preorder
$S-$extension of the partial order $\preceq.$

Moreover, for any pair of points $a,\,b \in X$ such that $(a,\,b)\not \in \preceq\cup S$
and $(b,\,a) \not \in \preceq\cup S,$ there exists a total preorder $\precsim$ which is
a total preorder $S-$extension of the partial order $\preceq$ and $(a,\,b) \in
\precsim.$

The intersection of all total preorder $S-$extensions of a preorder $\preceq$ coincides
with the preorder $\preceq\cup S.$ }

\vspace{2mm}

Along with each relation of preorder $\precsim$ we will consider {\it the indifference
relation} $I_\precsim :=\{(x,\,y) \in X \times X\,|\, (x,\,y)\not \in P_\precsim,
(y,\,x) \not \in P_\precsim \}$ corresponding to $\precsim.$ In the general case the
indifference relation $I_{\precsim}$ is reflexive and symmetric, i.~e. $I_{\precsim}$ is
a tolerance relation. An indifference relation $I_{\precsim}$  is in addition transitive
(and, consequently, it is an equivalency relation in this case) if and only if
$P_\precsim$ is negatively transitive (it means that the negation of $P_\precsim,$ i.~e.
the relation $(X\times X) \setminus P_\precsim,$ is transitive). It immediately follows
from the definition of $I_\precsim$ that every pair of points $x,\,y \in X$ satisfies
one and only one of the following three alternatives: $(x,\,y) \in P_\precsim,\,\,(y,\,x)
\in P_\precsim$ and $(x,\,y) \in I_\precsim.$

Another binary relation on $X$ generated by a preorder $\precsim$ is {\it the
equipotency relation} $R_\precsim,$ which is defined by
$$
(x,\,y) \in R_\precsim \Longleftrightarrow \{z \in X\,|\,(x,\,z) \in I_\precsim\} = \{z
\in X\,|\,(y,\,z) \in I_\precsim\}.
$$
It is not hard to verify that $R_\precsim$ is an equivalency relation on $X$ with
$R_\precsim \subset I_\precsim.$  The equipotency relation $R_\precsim$ is equal to the
indifference relation $I_\precsim,$ i.~e. $R_\precsim = I_\precsim,$ if and only if
$P_\precsim$ is negatively transitive (or, equivalently, if and only if $I_\precsim$ is
transitive).

\vspace{3mm}

{\bf Proposition 1}\,\, {\it  Let $\precsim$ be a preorder on $X.$ An equivalency
relation $S \subset X \times X$ holds the equalities $P_\precsim \circ S = S \circ
P_\precsim = P_\precsim$ if and only if  $S \subset R_\precsim.$}

\vspace{3mm}

\textit{Proof}\,\, Assume that an equivalency relation $S$ satisfies
the equalities $P_\precsim \circ S = S \circ P_\precsim = P_\precsim$ and let $(x,\,y)
\in S.$ The alternative $(x,\,y) \in P_\precsim$ is impossible, because otherwise it
would follow from $(y,\,x) \in S$ and from the equality $P_\precsim \circ S = P_\preceq$
that $(x,\,x) \in P_\precsim,$ but it contradicts the asymmetric property of
$P_\precsim.$ Similarly we can show that the alternative $(y,\,x) \in P_\precsim$ is
also impossible. Hence, $(x,y) \in I_\precsim.$

Let us prove that in fact $(x,y) \in R_\precsim.$ Choose an arbitrary element $z \in X$
such that $(x,\,z) \in I_\precsim$ and consider the ordered pair $(y,\,z) \in X \times
X.$ The alternatives $(y,\,z) \in P_\precsim$ and $(z,\,y) \in P_\precsim$ are
impossible, because otherwise it would follow from $(x,\,y) \in S$ and $P_\precsim \circ
S = S \circ P_\precsim = P_\precsim$ that $(x,\,z) \in P_\precsim,$ which contradicts
the choice of $z.$ Hence, $(y,\,z) \in I_\precsim$ and, consequently, $\{z \in
X\,|\,(x,\,z) \in I_\precsim\} \subset \{z \in X\,|\,(y,\,z) \in I_\precsim\}.$ The
converse inclusion is proved in the similar way. Thus, $\{z \in X\,|\,(x,\,z) \in
I_\precsim\} = \{z \in X\,|\,(y,\,z) \in I_\precsim\}$ and we conclude from the
definition of $R_\precsim$ that $(x,\,y) \in R_\precsim.$

To prove the converse statement we note that the inclusions
$P_\precsim \subset P_\precsim \circ S$ and $P_\precsim \subset S \circ P_\precsim$ follow from the reflexivity of the relation $S.$ So
we need to prove the converse inclusions. Let $(x,\,y) \in S \circ P_\precsim.$ Then
there exists an element  $z \in X$ such that $(x,\,z) \in S$ and $(z,\,y) \in
P_\precsim.$ Assume that $(y,\,x) \in P_\precsim.$ Due to the transitivity of
$P_\precsim$ we conclude from $(z,\,y) \in P_\precsim$ that $(z,\,x) \in P_\precsim,$
which contradicts $(x,\,z) \in S \subset R_\precsim\subset I_\precsim.$ Consequently,
$(y,\,x) \not \in P_\precsim.$ The assumption $(x,\,y) \in I_\precsim$ also leads to a
contradiction. Indeed, for $(x,\,z) \in S \subset R_\precsim$ we have due to the
definition of $R_\precsim$ that $(x,\,y) \in I_\precsim$ implies $(z,\,y)\in
I_\precsim,$ which contradicts $(z,\,y)\in P_\precsim.$ Hence, $(x,\,y) \not \in
I_\precsim$ and, consequently, the alternative $(x,\,y) \in P_\precsim$ is uniquely
possible. Thus, $S \circ P_\precsim \subset P_\precsim.$

The inclusion $P_\precsim \circ S \subset P_\precsim$ is proved in the similar way.
\hfill $\Box$

\vspace{3mm}

{\bf Corollary}\,\, {\it For every partial order $\preceq$ on the set $X$ and every
equivalency  relation $S$ on the same set $X$ such that $S \subset R_\preceq$ there
exists a total preorder $S-$extension of $\preceq.$}

\vspace{2mm}

Let us consider now the general case, that is the case when a partial order $\preceq$
and an equivalency relation $S$ do not necessarily satisfy equalities \eqref{e1}. We
begin with the following (evident) necessary condition for the existence of a total
preorder $S-$extension of a partial order $\preceq.$

\vspace{3mm}

{\bf Theorem 4}\,\, {\it Let $S$ be an equivalency relation on a
set $X.$ If for a partial order $\preceq \subset X \times X$ there exists a total
preorder $S-$extension then $S \subset I_\preceq.$}

\vspace{3mm}

\textit{Proof}\,\, Let $(x,\,y) \in S.$ If $(x,\,y) \in P_\preceq$ or $(y,\,x) \in
P_\preceq,$ then for any total preorder $S-$extension $\precsim$ of a partial order
$\preceq$ we would have $(x,\,y) \in P_\precsim$ or $(y,\,x) \in P_\precsim,$
respectively. However, since $P_\precsim \cap S = \varnothing,$ the both alternatives
are impossible and, hence, $(x,\,y) \in I_\preceq.$ \hfill $\Box$

\vspace{3mm}

{\bf Proposition 2}\,\, {\it Let $S$ be an equivalency relation on a set $X$ and
$\preceq$ a partial order defined on the same set $X.$ Then $S
\subset I_\preceq$ if and only if the composition $S \circ P_\preceq \circ S$ is
irreflexive.}

\vspace{3mm}

\textit{Proof}\,\, Recall that the irreflexivity of the binary relation $S \circ P_\preceq
\circ S$ means that $(x,\,x) \not \in S \circ P_\preceq \circ S$ for all $x \in X.$

Let $S \subset I_\preceq.$ Assume that contrary to the assertion of the
proposition the composition $S \circ P_\preceq \circ S$  is not irreflexive. The latter
means that $(x,\,x) \in S \circ P_\preceq \circ S$ for some $x \in X.$ Due to the
definition of the composition we can find such elements $y,z\in X$ that $(x,\,y) \in
S,\,\,(y,\,z) \in P_\preceq$ and $(z,\,x)\in S.$ Since $S$ is transitive, it follows
from $(z,\,x)\in S$ and $(x,\,y)\in S$ that $(z,\,y) \in S.$ Hence, since $S$ is
symmetric, $(y,\,z) \in S \cap P_\preceq,$ which contradicts $I_\preceq \cap
P_\preceq = \varnothing.$ This proves that $S \subset I_\preceq.$

Assume now that the composition $S \circ P_\preceq \circ S$ is
irreflexive, but the inclusion $S\subset I_{\preceq}$ is not the case. Then there exists
$(x,\,y) \in S$ such that $(x,y)\not\in I_{\preceq}$ and, consequently, either $(x,\,y)
\in P_\preceq,$ or $(y,\,x) \in P_\preceq.$ If $(x,y)\in P_\preceq,$ it follows from
$(y,\,x) \in S,\,\,(x,\,y) \in P_\preceq$ and $(y,\,y) \in S$ that $(y,\,y) \in S\circ
P_\preceq\circ S,$ but it is impossible since $S \circ P_\preceq \circ S$ is
irreflexive. Using the similar argument, we conclude that the case $(y,\,x) \in
P_\preceq$ is also impossible. It proves that $S$ is a subset of $I_{\preceq}.$
\hfill $\Box$

\vspace{3mm}

Recall that a binary relation $G \subset X\times X$ is said to be {\it acyclic} if for
any finite collection of elements $x_1,\,x_2,\ldots,x_n \in X$ it follows from
$(x_i,\,x_{i+1 }) \in G,\,i=1,\ldots,n-1,$ that $(x_n,\,x_1) \not \in G.$  {\it The
transitive hull} of a binary relation $G$ is the smallest transitive relation $TH(G)$
containing $G.$  There holds the equality $TH(G)=\cup\{G^n\,|\,n \in {\mathbb{N}}\},$
where ${\mathbb{N}}$ stands for the set of natural numbers and $G^n :=
\underbrace{G\circ G \circ \ldots \circ G}_n.$ It immediately follows from the latter
equality  that a binary relation $G$ is acyclic if and only if its transitive hull
$TH(G)$ is irreflexive (or, equivalently, if and only if $TH(G)$ is asymmetric).

\vspace{3mm}

{\bf Theorem 5}\,\, {\it Let $S$ be an equivalency relation on a set $X$ and
$\precsim$ a partial order defined on the same set $X.$ Then the
following statements are equivalent:

$(i)$ there exists a total preorder $S-$extension of $\preceq;$

$(ii)$ the composition $S \circ P_\preceq \circ S$ is acyclic.}

\vspace{3mm}

\textit{Proof}\,\, $(i) \Longrightarrow (ii)$ Let a total preorder $\precsim$ be a total
preorder $S-$extension of a partial order $\preceq.$ Assume that the composition $S\circ
P_\preceq \circ S$ is not acyclic and let the collection $x_1,\,x_2,\ldots,x_m \in X,\,m
\ge 2,$ hold $(x_i,\,x_{i+1}) \in S\circ P_\preceq \circ S,\,i=1,\ldots,m-1,$ and
$(x_{m},\,x_1) \in S\circ P_\preceq \circ S.$ Then there exist collections
$y_1,\,y_2,\ldots,y_m \in X$ and $z_1,\,z_2,\ldots,z_m \in X$ such that $(x_i,\,y_i) \in
S,\,(y_i,\,z_i) \in P_\preceq,\,(z_i,\,x_{i+1}) \in S,\, i=1,2,\ldots,m-1,$ and
$(x_m,\,y_m) \in S,\,(y_m,\,z_m) \in P_\preceq,\,(z_m,\,x_{1}) \in S.$ The inclusion
$P_\preceq \subset P_\precsim$ implies that $(y_i,\,z_i) \in P_\precsim,\,i=1,\ldots,m.$
Since $S= R_\precsim,$ we conclude from Proposition 1 that $P_\precsim \circ S = S \circ
P_\precsim = P_\precsim.$ Hence, it follows from $(x_i,\,y_i) \in S,\,(y_i,\,z_i) \in
P_\precsim,\,(z_i,\,x_{i+1}) \in S,\, i=1,2,\ldots,m-1,$ that $(x_i,\,x_{i+1}) \in
P_\precsim,\, i=1,2,\ldots,m-1,$ whence, due to the transitivity of $P_\precsim,$ we
obtain $(x_1,\,x_m) \in P_\precsim.$ On the other hand, from $(x_m,\,y_m) \in
S,\,(y_m,\,z_m) \in P_\precsim,\,(z_m,\,x_{1}) \in S,$ using the equalities $P_\precsim
\circ S = S \circ P_\precsim = P_\precsim,$ we deduce $(x_m,\,x_1) \in P_\precsim.$ This
is a contradiction because $P_\precsim$ is asymmetric. It proves that $S \circ P_\preceq
\circ S$ should be acyclic.

$(ii) \Longrightarrow (i)$ The relation $S \circ P_\preceq \circ S$ is acyclic if and
only if its transitive hull $TH(S \circ P_\preceq \circ S)$ is asymmetric. Using the
equality $TH(S \circ P_\preceq \circ S) = \cup\{(S \circ P_\preceq \circ S)^n\,|\,n \in
{\mathbb{N}}\}$ it is not difficult to verify that $S \circ TH(S \circ P_\preceq \circ
S)= TH(S \circ P_\preceq \circ S)\circ S = TH(S \circ P_\preceq \circ S).$ Hence, due to
Theorem~3 there exists a total preorder $\precsim$ on $X$ which is a total preorder
$S-$extension of a partial order $E\cup TH(S \circ P_\preceq \circ S)$ (recall that
$E:=\{(x,\,x) \in X \times X\,|\, x \in X\}$ is the equality relation on $X$). Since
$P_\preceq \subset TH(S \circ P_\preceq \circ S),$ the preorder $\precsim$ is also a
total preorder $S-$extension of a partial order $\preceq.$ \hfill $\Box$

\vspace{3mm}

{\bf Theorem 6}\,\, {\it Let $\preceq$ be a partial order on a set $X$ and
$S$ an equivalency relation defined on the same set $X.$  If the
composition $S \circ P_\preceq \circ S$ is acyclic then the intersection of all total
preorder $S-$extension of the partial order $\preceq$ is the preorder $S \cup TH(S \circ
P_\preceq \circ S),$ that is the preorder whose asymmetric part is the transitive hull
of $S \circ P_\preceq \circ S$ and whose symmetric part coincides with $S.$}

\vspace{0mm}

\textit{Proof}\,\, Since the asymmetric part of every total preorder $S-$extension of
the partial order $\preceq$ is transitive and contains $S \circ P_\preceq \circ S,$ it also contains the
transitive hull of $S \circ P_\preceq \circ S.$ Hence, every total preorder
$S-$extension of the partial order $\preceq$ is at the same time a total preorder
$S-$extension of the partial order $E \cup TH(S \circ P_\preceq \circ S).$ Conversely,
it follows from $P_\preceq \subset S \circ P_\preceq \circ S$ that every total preorder
$S-$extension of the partial order $E \cup TH(S \circ P_\preceq \circ S)$ is a total
preorder $S-$extension of the partial order $\preceq.$ Since $S \circ TH(S \circ
P_\preceq \circ S)= TH(S \circ P_\preceq \circ S)\circ S = TH(S \circ P_\preceq \circ
S),$ we conclude from the second statement of  Theorem 3 that the intersection of all
total preorder $S-$extension of the partial order $\preceq$ coincides with the preorder
$S \cup TH(S \circ P_\preceq \circ S).$ \hfill $\Box$

\vspace{3mm}

Given a partial order $\preceq$ on $X,$ by the symbol $\Sigma(\preceq)$ (respectively,
$\Sigma^\ast(\preceq)$) we denote the collection consisting of all equivalency
relations $S$ defined on $X$ such that the composition $S\circ P_\preceq \circ S$ is
irreflexive (respectively, $S\circ P_\preceq \circ S$ is acyclic). Clearly,
$\Sigma^\ast(\preceq)$ is a subcollection of the collection $\Sigma(\preceq).$ It also
follows from Proposition 2 and Theorem 5 that $S \in \Sigma(\preceq)$ if and only if $S
\subset I_\preceq$ and $S \in \Sigma^\ast(\preceq)$ is equivalent to the existence of a
total preorder $S-$extension of the partial order $\preceq.$

\vspace{2mm}

{\bf Theorem 7} {\it Let  $S$ be an equivalency
relation defined on a set $X.$
A partial order $\preceq$ defined on the same set $X$ has a
unique total preorder $S-$extension if and only if
$S$ is maximal (in inclusion) in the subcollection $\Sigma^\ast(\preceq)$ and the transitive hull of the composition $S \circ P_\preceq \circ S$ is negatively transitive.}

\vspace{2mm}

\textit{Proof} Let $\precsim$ be a unique total preorder $S-$extension of the partial
order $\preceq\!.$ Suppose to the contrary that $S$ is not maximal (in inclusion) in the
subcollection $\Sigma^\ast(\preceq).$ Then there exists an equivalency relation $S'$ in
$\Sigma^\ast(\preceq)$ such that $S \subset S',\,S \ne S'.$ Let $\precsim'$ be an
arbitrary total preorder $S'$--extension of the partial order $\preceq.$ Suppose that
$P_{\precsim'} \not\subset P_\precsim.$ Denote by $\precsim^*$ an arbitrary total
preorder $S-$extension of the partial order $P_{\precsim'}\cup E.$ The existence of
$\precsim^*$ follows from the inclusion $S \subset S'=R_{{\precsim'}}$ and Corollary 6.
Since $P_\preceq \subset P_{\precsim'} \subset P_{\precsim^*},$ then $\precsim^*$ is
also a total preorder $S-$extension of the initial partial order $\preceq\!.$ It follows
from the assumption $P_{\precsim'} \not\subset P_\precsim$ and the inclusion
$P_{\precsim'} \subset P_{\precsim^*}$ that $P_{\precsim^*} \not\subset P_\precsim.$
Hence $\precsim^* \ne \precsim.$ Since it contradicts the uniqueness of a total preorder
$S-$extension of the partial order $\preceq,$ then the inclusion $P_{\precsim'}
\not\subset P_\precsim$ is impossible and, consequently, we have $P_{\precsim'} \subset
P_\precsim.$ In this case we define on $X$ the relation $\precsim^\circ:=
P_{\precsim'}\cup (P^{-1}_{\precsim}\cap S')\cup S.$ It is not difficult to verify that
$\precsim^\circ$ is a total preorder, which differs from $\precsim$ only on the
equivalency classes of $S',$ where it coincides with the converse relation of
$\precsim.$ Since $P_\preceq \subset P_\precsim' \subset P_{\precsim^\circ} :=
P_{\precsim'}\cup (P^{-1}_{\precsim}\cap S')$ and $S_{\precsim^\circ}= S,$ then
$\precsim^\circ$ is a total preorder $S-$extension of the partial order $\preceq$ which
differs from $\precsim'.$ Again we get the contradiction to the uniqueness of a total
preorder $S-$extension of the partial order $\preceq.$ This completes the proof that $S$
is maximal (in inclusion) in the subcollection $\Sigma^\ast(\preceq).$

It remains to prove that $TH(S \circ P_\preceq \circ S)$ is negatively transitive.
Since $\precsim$ is the unique total preorder $S-$extension of the partial order
$\preceq$ then due to Theorem 6 we conclude that $P_\precsim = TH(S \circ P_\preceq \circ S).$
Notice now that $P_\precsim$ is the asymmetric part of the total preorder $\precsim$ and
therefore it is negatively transitive. Hence, $TH(S \circ P_\preceq \circ S)$ is negatively
transitive too.

For the converse, notice that the assumption $S \in \Sigma^\ast(\preceq)$ is equivalent
to the asymmetry property of the transitive hull $TH(S\circ P_\preceq \circ S).$ It
implies that the relation $\precsim:=S\cup TH(S\circ P_\preceq \circ S)$ is the
preorder. Since  $TH(S\circ P_\preceq \circ S)$ is negatively transitive, the
indifference relation $I_\precsim$ corresponding to $\precsim$ is transitive and hence
$I_\precsim$ is an equivalency relation. Then the relation  $I_\precsim \cup TH(S\circ
P_\preceq \circ S)$ is a total preorder and, moreover, it follows from $P_\preceq
\subset TH(S\circ P_\preceq \circ S$ that $I_\precsim \in \Sigma^\ast(\preceq).$ Thus,
the preorder $\precsim:=S\cup TH(S\circ P_\preceq \circ S)$ is total and consequently it
is a total preorder $S-$extension of the partial order $\preceq.$ From Theorem 6 we
conclude that there are no other total preorder $S-$extensions of the partial order
$\preceq.$ \hfill $\Box$

\vspace{1mm}

{\bf Theorem 8}\,\, {\it Let $\preceq$ be a partial order on $X$ and $S$ an equivalency
relation defined on the same set $X.$
The relation $S\cup\!(S\circ P_\preceq \circ S)$ is the unique
total preorder $S-$extension of the partial order $\preceq$ if and only if
$S$ belongs to the subcollection
$\Sigma^\ast(\preceq)$ and is maximal (by inclusion) in the collection
$\Sigma(\preceq).$ }

\vspace{3mm}

\textit{Proof}\,\,Assume that an equivalence relation $S$ belongs to $\Sigma^\ast(\preceq)$ and is maximal (by inclusion) in $\Sigma(\preceq).$ First we prove that for any $x,\,y \in X$ there holds exactly one alternative
of the following three ones:
$$(x,\,y) \in S,\,\,(x,\,y) \in S\circ P_\preceq \circ
S,\,\,(y,\,x) \in S\circ P_\preceq \circ S.$$

Choose $x,\,y \in X$ with $(x,\,y) \not \in S.$  If $(x,\,y) \in P_\preceq$ or $(y,\,x)
\in P_\preceq,$ then $(x,\,y) \in S\circ P_\preceq \circ S$ or $(y,\,x) \in S\circ
P_\preceq \circ S,$ respectively. Let $(x,\,y) \in I_\preceq \setminus S$ and let
$[x]_S$ and $[y]_S$ be the equivalency classes of $S$ containing $x$ and $y,$
respectively. Since $S$ is maximal (in inclusion) in $\Sigma(\preceq),$ there exist $x_1
\in [x]_S$ and $y_1 \in [y]_S$ such that either $(x_1,\,y_1) \in P_\preceq,$ or
$(y_1,\,x_1) \in P_\preceq.$ Indeed, if $(x_1,\,y_1) \not \in P_\preceq$ and
$(y_1,\,x_1) \not \in P_\preceq$ for any $x_1 \in [x]_S$ and $y_1 \in [y]_S,$ then the
equivalency relation $S'$ defined on $X$ by
$$
(u,v) \in S' \Longleftrightarrow \,\,\text{either}\,\,(u,\,v) \in
S\,\,\text{or}\,\,u,\,v \in [x]_S\cup [y]_S.
$$
satisfies $S' \subset I_\preceq.$  Since $S \subset S',\,S \ne S',$ it contradicts
maximality (by inclusion) of $S$ in $\Sigma(\preceq).$

Thus, for any pair $(x,\,y) \in I_\preceq \setminus S$ there exist $x_1 \in [x]_S$ and
$y_1 \in [y]_S$ such that either $(x_1,\,y_1) \in P_\preceq,$ or $(y_1,\,x_1) \in
P_\preceq.$

If $(x_1,\,y_1) \in P_\preceq$ is the case then it follows from $(x,\,x_1) \in
S,\,\,(x_1,\,y_1) \in P_\preceq,\,\,(y_1,\,y) \in S$ that $(x,\,y) \in S\circ P_\preceq
\circ S.$ Similarly, in the case when $(y_1,\,x_1) \in P_\preceq$ we get from $(y,\,y_1)
\in S,\,\,(y_1,\,x_1) \in P_\preceq,\,\,(x_1,\,x) \in S$ that $(y,\,x) \in S\circ
P_\preceq \circ S.$

The fact that for any $x,\,y \in X$ there holds exactly one of the three possible
alternatives follows from the assumption $S \in \Sigma^\ast(\preceq)$ or, equivalently
from the acyclicity of $S\circ P_\preceq \circ S.$

Assume now that $\precsim\, \subset\, X \times X$ is an arbitrary
total preorder $S-$extension of the partial order $\preceq$ (the existence of total
preorder $S-$extensions for $\preceq$ is guaranteed by the assumption that $S \in
\Sigma^\ast(\preceq).$) It follows from $S_\precsim = S$ и $P_\preceq \subset
P_\precsim$ that $S \circ P_\preceq \subset S \circ P_\precsim = P_\precsim$ and then $S
\circ P_\preceq \circ S \subset P_\precsim \circ S = P_\precsim.$ To prove the converse
inclusion let us consider a pair $(x,\,y) \in P_\precsim.$ Then $(x,\,y) \not \in S$ and
by the assertion proved above we have either $(x,\,y) \in S \circ P_\preceq \circ S,$ or
$(y,\,x) \in S \circ P_\preceq \circ S.$ The latter is impossible because it contradicts
the asymmetry property of $P_\precsim.$ Hence, $(x,\,y) \in S \circ P_\preceq \circ S$
and we get $P_\precsim = S \circ P_\preceq \circ S.$ Since $\precsim$ is an arbitrary
total preorder $S-$extension of $\preceq$ we conclude that $S \cup (S \circ P_\preceq
\circ S)$ is the unique total preorder $S-$extension of the partial preorder $\preceq.$

To verify the converse, assume that $S\cup\!(S\circ P_\preceq \circ S)$ is a total preorder $S-$extension of a partial order $\preceq\!.$ Obviously, $S \in \Sigma^\ast(\preceq).$ Suppose that $S \subset S'$ for some $S' \in \Sigma(\preceq).$ Then $S\circ P_\preceq \circ S \subset S'\circ P_\preceq \circ S'$ and $S' \cap (S'\circ P_\preceq \circ S') = \varnothing.$ Since the preorder $S\cup\!(S\circ P_\preceq \circ S)$ is total, we get that $S' \subset S.$ Hence, $S = S'$ and it proves that $S$ is maximal in $\Sigma(\preceq).$
\hfill $\Box$

\vspace{2mm}

{\bf Remark}\,\, Let $X$ be a real vector space and let a partial order $\preceq$ and an
equivalency relation $S$ defined on $X$ be compatible with algebraic operations on $X.$
Then the condition $S \cap P_{\preceq} = \varnothing$ is both necessary and sufficient
for the existence of a compatible total preorder $S-$extension of $\preceq.$  This
criterion follows from the Kakutani--Tukey theorem on separation of convex sets by
halfspaces (see, for instance, [12, Theorem 1.9.1, p. 12]) and from the duality between
compatible total preorders and conic halfspaces [15--16].


\vspace{5mm}

{\small

\textbf{References}

\vspace{5mm}

1. Szpilrajn, E.: Sur l'extension de l'ordre partial. Fundamenta Mathematicae.
{\bf 1}, 386--389 (1930)

2. Dushnik, B., Miller, E.: Partially ordered sets. American Journal of
Mathematics. {\bf 63}, 600--610 (1941)

3. Fishbern, P.C.: Utility Theory for Decision Making. John Wiley
\& Sons, Inc., New York (1970)

4. Harzheim, E.: Odered Sets. Springer, New York (2005)

5. Fuchs, L.: Partially Ordered Algebraic Systems. Pergamon Press, Oxford (1963)

6. Klee, V.: The structure of  semispaces. Math. Scand. {\bf 4}, 54--64 (1956)

7. Bosi, G., Herden, G.: On a strong continuous analogue of the Szpilrajn theorem
and its strengthening by Dushnik and Miller, Order. {\bf 22}, 329--342 (2005)

8. Bosi, G., Herden, G.: On a possible continuous analogue of the Szpilrajn
theorem and its strengthening by Dushnik and Miller. {\bf 23}, 271--296 (2006)

9. Debreu, G.: Representation of preference ordering by a numerical function.
In: Thrall, R., Coombs, C.C., Davis, R. (eds.) Decision Processes. pp.\, 159~--~166. Wiley, New York (1954)

10. Kiruta, A. Ya., Rubinov, A.M., Yanovskaja, E.B.: Optimalnyj vybor raspredelenij v socialnykh i ekonomicheskikh sistemakh (Optimal choice of distribution
in social and economic systems).--  Nauka, Leningrad Department, Leningrad (1980)

11. Bridges, D.S., Mehta, G.B.: Representations of Preferences Orderings. Lecture
Notes in Economics and Mathematical Systems, vol. 422.  Springer Verlag, Berlin (1995)

12. Hille, E., Phillips, R.S.: Functional analysis and semi--groups. AMS
Colloquium Publications, vol. 31, Providence, R. I. (1957)

13. Martinez-Legaz, J.E., Singer, I.: Lexicographical Separation  in ${\mathbb{R}}^{n}.$ Linear  Algebra  and   Appl. {\bf 90}, 147--163 (1987)

14. Martinez-Legaz,  J.E.,  Singer,  I.: Compatible  Preorders  and Linear Operators
on ${\mathbb{R}}^{n},$  Linear  Algebra  and Appl. {\bf 153}, 53--66 (1991)

15. Gorokhovik, V.V., Shinkevich, E.A.: Teoremy ob otdelimosti vypuklykh mnozhestv
stupenchato--linejnymi funkcijami i ikh prilozhenija k zadacham vypukloj opimizacii
(Theorems on separation of convex sets by step-linear functions and their applications
to convex optimization problems). In: Gorokhovik, V.V., Zabreiko, P. P. (eds.)
Nonlinear Analysis and Applications, Trudy Instituta Matematiki NAN Belarusi.
vol.~1, pp. 58--85 (1998) 

16. Gorokhovik, V.V., Shinkevich, E.A.: Geometric structure and classification of
infinite-dimensional halfspaces.  In: Przeworska--Rolewicz, D. (ed.) Algebraic Analysis and Related Topics, Banach Center Publica\-tions, vol.~53, pp. 121--138. Institute of Mathematics PAN, Warsaw (2000)

}

\end{document}